\documentclass[fleqn]{mat01}
\usepackage{times,mathtimy,amssymb,epsfig}
\begin{document}

\setcounter{page}{443}
\firstpage{443}

\def\fr{\frac}
\def\p{\partial}
\def\u{{\bf u}}
\def\d{\delta}
\def\e{\epsilon}
\def\bd{{\bf {\cdot}}}

\def\thoe{\trivlist\item[\hskip\labelsep{{\bf Theorem}}]}
\def\thoee{\trivlist\item[\hskip\labelsep{{\bf Theorem.}}]}
\newtheorem{theo}{\bf Theorem}
\renewcommand\thetheo{\arabic{theo}}
\def\remar{\trivlist\item[\hskip\labelsep{{\it Remark.}}]}

\markboth{Ganesh C Gorain and Sujit K Bose}{Uniform stability of damped nonlinear vibrations}

\title{Uniform stability of damped nonlinear vibrations\\ of an elastic
string}

\author{GANESH C GORAIN and SUJIT K BOSE$^{*}$}

\address{Department of Mathematics, J.K. College, Purulia~723~101, India\\
\noindent $^{*}$BE 188, Salt Lake City, Kolkata~700~064, India\\
\noindent Email: sujitkbose@yahoo.com}

\volume{113}

\mon{November}

\parts{4}

\Date{MS received 21 April 2003; revised 5 August 2003}

\begin{abstract}
Here we are concerned about uniform stability of damped nonlinear
transverse vibrations of an elastic string fixed at its two ends. The
vibrations governed by nonlinear integro-differential equation of
Kirchoff type, is shown to possess energy uniformly bounded by
exponentially decaying function of time. The result is achieved by
considering an energy-like Lyapunov functional for the system.
\end{abstract}


\maketitle

\section{Introduction}

The textbook treatment of transverse vibrations of a metallic wire
(string), governed by the linear wave equation, does not portray
existence of whirling out of plane motion. The phenomenon was first
observed by Hunton as reported by Harrison \cite{9}. The reference to
other subsequent experiments can be found in the book by Nayfeh and Mook
\cite{14}. It is observed that the whirling motion, sometimes also
referred to as {\em ballooning motion} occurs when the amplitude and
frequency of a plane excitation and the phase difference between the
response and the excitation exceeds certain critical values. The
explanation of the phenomenon lies in a nonlinear treatment of the
problem. Long ago Kirchoff \cite{10}, taking longitudinal elasticity
also into consideration, gave the following nonlinear
integro-differential equation for transverse vibrations confined to a
plane:
\begin{subequations}
\begin{equation}
\fr{\p^2 u}{\p t^2} + 2\d\fr{\p u}{\p t} = \left[\,a^2 +
b\int_0^l\left|\fr{\p u}{\p x}\right|^2 {\rm d}x\right]\,\fr{\p^2 u}{\p
x^2},\quad 0 < x < l,\quad t \ge 0
\end{equation}
where $u(x, t)$ is the transverse deflection, $l$ the length of the
string, $\d > 0$ the coefficient of damping, $a^2 = T_0/\rho A$ and $b =
E/2\rho l$. The constants $T_0, \rho, A, E$ are respectively the initial
axial tension, the mass density, the cross-section area and Young's
modulus of the string. For treating nonlinear out of plane motion with
components $(v, w)$ in the perpendicular $y$ and $z$ directions, Anand
\cite{1} derives a generalization of equation (1a):
\begin{equation}
\fr{\p^2\u}{\p t^2} + 2\d\fr{\p\u}{\p t} = \left[\,a^2 +
b\int_0^l\left|\fr{\p\u}{\p x}\right|^2 {\rm d}x\right]\,\fr{\p^2\u}{\p
x^2},\quad 0 < x < l,\quad t \ge 0
\end{equation}
\end{subequations}
where $\u$ has components $(v, w)$. A systematic rederivation of (1b)
can also be found in Nayfeh and Mook \cite{14}. The book \cite{14} also
refers to other authors who derived different nonlinear differential
equations under various assumptions.

For a string  fixed at both ends, we have
\begin{equation}
\u(0, t) = \u(l, t) = 0
\end{equation}
and if $\u_0(x),\;\u_1(x)$ are initial displacement and velocity
\begin{equation}
\u(x, 0) = \u_0(x),\quad\fr{\p\u}{\p t}(x,0) = \u_1(x),\quad 0 < x < l.
\end{equation}
The existence and uniqueness in the scalar case are discussed in
\cite{12} and \cite{15}.

Arosio and Spagnola \cite{3}, Gough \cite{7} and several other authors
referred to in Nayfeh and Mook \cite{14} have studied the problem (1a, b)
from mathematical and physical points of view with a focus on modal
analysis and eigenvalues of the system. In particular, Anand \cite{2}
has investigated (1b) proving that the vibrations are stable even in the
planar case. The elaborate method of treatment employs stability
analysis of the time-dependent Hill's equation satisfied by the modal
amplitudes. In recent years, Shahruz \cite{16} has treated the problem
(1a) with a bounded input disturbance and his result shows boundedness
of the ouput. The approach is systems and control theoretic.

The mathematical theory of stabilisation of distributed parameter
systems is currently of interest in view of application to vibration
control of various structural elements. The most common class of
vibration control mechanism is of passive type that absorbs vibration
energy. A system is called strongly stable if the total energy $E(t)$ of
each solution of the system converges to zero as time $t\rightarrow
+\,\infty$. If the convergence is uniform for $t > 0$ with respect to
all initial data in the energy space for which $E(0) < \infty$, the
system is called uniformly stable. Here our investigation is a direct
Lyapunov stability approach to obtain uniform exponential energy decay
estimate for the system (1b), (2), (3) and establish that the system is
a passive energy absorber. Such estimate has earlier been obtained by
Gorain \cite{4} for damped linear wave equation in a bounded domain in
$R^n$. Gorain and Bose \cite{5,6} have obtained such estimates for an
internally damped beam governed by linear equations in torsional and
flexural modes of vibration.

\section{Energy of the system}

The total energy at time $t$, of the system (1b)--(3) is defined by
the functional
\begin{equation}
E(t) = \fr{1}{2}\int_0^l\left[\left|\fr{\p\u}{\p t}\right|^2 +
a^2\left|\fr{\p\u}{\p x}\right|^2\right]{\rm d}x +
\fr{b}{4}\left[\int_0^l\left|\fr{\p\u}{\p x}\right|^2{\rm
d}x\right]^2,\quad t \ge 0.
\end{equation}
Diferentiating (4) with respect to $t$ and using the governing equation
(1b) we obtain
\begin{align*}
\fr{{\rm d}E}{{\rm d}t} &= \int_0^l\fr{\p\u}{\p t}\,\bd\,\left[\left(a^2
+ b\int_0^l\left|\fr{\p\u}{\p x}\right|^2{\rm d}x\right)\,\fr{\p^2\u}{\p
x^2} - 2\d\fr{\p\u}{\p t}\right]{\rm d}x\\[.2pc]
&\quad\ + \left(a^2 + b\int_0^l\left|\fr{\p\u}{\p
x}\right|^2{\rm d}x\right)\,\int_0^l\fr{\p\u}{\p x}\,\bd\,\fr{\p^2\u}{\p x\p
t}\;{\rm d}x\\[.2pc]
&= \left(a^2 + b\int_0^l\left|\fr{\p\u}{\p
x}\right|^2{\rm d}x\right)\int_0^l\fr{\p}{\p x}\left(\fr{\p\u}{\p
t}\,\bd\,\fr{\p\u}{\p x}\right){\rm d}x - 2\d\int_0^l\left|\fr{\p\u}{\p
t}\right|^2{\rm d}x.
\end{align*}
The value of the first integral vanishes in view of the boundary
conditions (2). We thus get
\begin{equation}
\fr{{\rm d}E}{{\rm d}t} = - 2\d\int_0^l\left|\fr{\p\u}{\p
t}\right|^2{\rm d}x,\quad t \ge 0.
\end{equation}
The negativity of the right hand side of (5) shows that the energy
$E(t)$ is decreasing with time due to the incorporation of passive viscous
damping coefficient $\d > 0$. The system is thus non-energy conserving.
Naturally the question arises as to whether the system is uniformly
stable or not. An affirmative answer is contained in the following
theorem, which essentially states that in fact the energy of vibration
has uniform exponential decay.

\begin{thoee} {\it Let $\u(x, t)$ be a solution of the system {\rm (1b)--(3)}
with $\u_0\in H^1[\,0, l]$ and $\u_1\in L^2[\,0, l]${\rm ,} where $H^1[\,0,
l]$ is the Sobolev space of order {\rm 1,} then the energy $E(t)$ of the
system defined by {\rm (4)} satisfies
\begin{equation}
E(t) \le M {\rm e}^{-\mu t} E(0),\quad t > 0
\end{equation}
for some reals $\mu > 0$ and $M > 1$.}\vspace{.5pc}
\end{thoee}

From eq.~{\rm (4)}
\begin{equation}
E(0)=\frac{1}{2}\int_0^l\left[|\u_1|^2+a^2\left|\frac{{\rm d}\u_0}{{\rm
d}x}\right|^2\right]{\rm d}x + \frac{b}{4}\left[\int_0^l\left|\frac{{\rm
d}\u_0}{{\rm d}x}\right|^2{\rm d}x\right]^2.
\end{equation}
From eq.~(6) it follows that if $E(0)<\infty$ with respect to all
initial values $\u_0\in H^1[0,l]$ and $\u_1\in L^2[0,l]$ then
$E(t)\rightarrow 0$ as $t\rightarrow \infty$ and the system is uniformly
stable. In order to prove (6), let $\e > 0$ be a fixed constant.
Proceeding as in Gorain \cite{4} (see also Komornik \cite{11}), we define an
energy-like Lyapunov functional $V$ according to
\begin{equation}
V(t) = E(t) + \e G(t),
\end{equation}
where
\begin{equation}
G(t) = \int_0^l\left(\u\,\bd\,\fr{\p\u}{\p t} +
\d|\u|^2\right){\rm d}x.
\end{equation}
Now since $\u(0, t) = {\bf 0} = \u(l, t)$, the Poincar\'{e}-like
Scheefer's inequality
\begin{equation}
\int_0^l|\u|^2{\rm d}x \le \fr{l^2}{\pi^2}\int_0^l\left|\fr{\p\u}{\p
x}\right|^2{\rm d}x,\quad t > 0
\end{equation}
holds (see \cite{13} and \cite{8}). A simple proof is provided in the
Appendix. Hence using Schwarz's inequality, we obtain
\begin{align}
\left|\int_0^l\u\,\bd\,\fr{\p\u}{\p t}{\rm d}x\right| &\le \fr{l}{\pi
a}\int_0^l\left|\fr{\p\u}{\p t}\right|\left|\fr{\pi a}{l}\u\right|{\rm
d}x\nonumber\\[.2pc]
&\le \fr{1}{2}\fr{l}{\pi a}\int_0^l\left[\left|\fr{\p\u}{\p
t}\right|^2 + \fr{\pi^2 a^2}{l^2}|\u|^2\right]{\rm d}x\nonumber\\[.2pc]
&\le \fr{l}{\pi a}E(t),\quad t > 0
\end{align}
where (10) and the defining equation (4) are used. Also, similarly
\begin{equation}
0 \le \d\int_0^l|\u|^2{\rm d}x \le \fr{\d l^2}{\pi^2
a^2}\int_0^la^2\left|\fr{\p\u}{\p x}\right|^2{\rm d}x \le \fr{2\d l^2}{\pi^2
a^2}E(t),\quad t > 0.
\end{equation}
Thus using (11) and (12) in eq.~(9)
\begin{equation}
|G(t)| \le \fr{l}{\pi a}\left(1 + \fr{2\d l}{\pi a}\right)E(t) =
\mu_0E(t)\quad({\rm say}),\quad t > 0.
\end{equation}
Also from (9),
\begin{equation}
G(t) \ge \int_0^l\u\,\bd\,\fr{\p\u}{\p t}{\rm d}x \ge
-\left|\int_0^l\u\,\bd\,\fr{\p\u}{\p t}{\rm d}x\right| \ge -\, \fr{l}{\pi
a}E(t),\quad t > 0.
\end{equation}
The inequalities (13) and (14) yield for $V(t)$ (defined by (8)) the
estimates
\begin{equation}
\left(1 - \fr{\e\, l}{\pi a}\right)E(t) \le V(t) \le (1 +
\e\mu_0)E(t),\quad t > 0
\end{equation}
where we assume that $\e < \pi a/l$, so that the left hand side of (15)
is positive.

Next, differentiating (9) with respect to $t$, and using the governing
equation (1b), we obtain
\begin{equation*}
\fr{{\rm d}G}{{\rm d}t} = \int_0^l\left[\u\,\bd\,\fr{\p^2\u}{\p x^2}\left(a^2 +
b\int_0^l\left|\fr{\p\u}{\p x}\right|^2{\rm d}x\right) +
\left|\fr{\p\u}{\p t}\right|^2\right]{\rm d}x.
\end{equation*}
Integration by parts of the first term yields
\begin{align}
\fr{{\rm d}G}{{\rm d}t} &= -\,\left(a^2 + b\int_0^l\left|\fr{\p\u}{\p
x}\right|^2{\rm d}x\right)\int_0^l\left|\fr{\p\u}{\p x}\right|^2{\rm d}x +
\int_0^l\left|\fr{\p\u}{\p t}\right|^2{\rm d}x\nonumber\\[.2pc]
&\le -\,2\,E(t) + 2\int_0^l\left|\fr{\p\u}{\p t}\right|^2{\rm d}x,\quad
t > 0.
\end{align}
If we now differentiate (8) with respect to $t$ and then insert the
results of (5) and (16), we get
\begin{equation*}
\fr{{\rm d}V}{{\rm d}t} \le 2\,(\e - \d)\int_0^l\left|\fr{\p\u}{\p
t}\right|^2{\rm d}x - 2\,\e\,E(t),\quad t > 0
\end{equation*}
and so if we choose $\e < \mbox{min}\{\,\d, \pi a/l\}$,
\begin{equation*}
\fr{{\rm d}V}{{\rm d}t} \le -\,2\,\e\,E(t) \le -\,\fr{2\,\e}{1 +
\e\mu_0}\;V(t),\quad t > 0
\end{equation*}
where (15) has been used. Defining
\begin{equation}
\mu = \fr{2\,\e}{1 + \e\mu_0},
\end{equation}
the above inequality becomes
\begin{equation*}
\fr{{\rm d}V}{{\rm d}t} + \mu\,V(t) \le 0,\quad t > 0
\end{equation*}
and so multiplying by $\exp(\mu t)$ and integrating from 0 to $t$, we
obtain
\begin{equation*}
V(t) \le {\rm e}^{-\,\mu t} V(0),\quad t > 0.
\end{equation*}
Finally, using (15) again, we obtain the estimate (6) where
\begin{equation}
M = \fr{1 + \e\,\mu_0}{1 - \e\,l/\pi a} > 1.
\end{equation}
Hence the theorem.

\begin{remar} In the above we have also obtained exponential energy
decay estimate for the solution of the nonlinear damped elastic string
for initial conditions $\u_0\in H^1[\,0, l],\;\u_1\in L^2[\,0, l]$. The
decay rate $\mu$ is a function of $\e$ as given in eq.~(17), where $\e$
is the lesser of the two quantities $\delta$ and $\pi a/l$. Since ${\rm
d}\mu/{\rm d}\e = 2/(1 + \e\mu_0)^2 > 0, \;\mu $ is maximum when $\e$ is
maximum. In actual cases $\d$ is the smaller parameter unless the string
is very long. So from the restriction on $\e$, we can assume $\e \le \d
\le \pi a /l$. Hence $\e_{\max} = \d$ and so
\begin{equation}
\mu_{\max} = \fr{2\,\d}{1 + \d\mu_0} = \fr{2\,\d}{1 + \fr{\d l}{\pi a}(1
+ \fr{2\d l}{\pi a})}.
\end{equation}
The right hand side is a function of $\d l/\pi a \le 1$ and is maximum
when $\d l/\pi a = 1/\sqrt{2}$. Hence
\begin{equation}
\mu_{\max} \le \fr{2}{1 + 2\sqrt{2}}\fr{\pi a}{l} = 0.52\,\fr{\pi a}{l}
\end{equation}
for all possible cases of damping.

It is of some interest to note that Anand \cite{1} in his modal solution of
the problem assumed $\mu=\d$. This assumption holds when $\d<\pi a/l$.
In the contrary case $\mu=\pi a/l$ will suffice exponential decay of the
solution.

The quantity $M$ given by eq.~(18) however also increases with $\e$,
because ${\rm d}M/{\rm d}\e = (\mu_0 + l/\pi a)/(1 - \e l/\pi a)^2 > 0$.
Thus for $\e_{\max} = \d$,
\begin{equation}
M = \fr{1 + \fr{\d l}{\pi a}(1 +\fr{2\d l}{\pi a})}{1 - \d l/\pi a} =
\fr{1 + 2\sqrt{2}}{\sqrt{2} - 1} = 9.3
\end{equation}
when $\d l/\pi a = 1/{\sqrt2}$. Thus if we look for steeply falling
estimate, the energy gets multiplied many times in the initial stages.
\end{remar}

\section{Quasi-steady amplitude estimate}

Scheefer's inequality (10) when applied to the energy functional $E(t)$
given by equation (4) yields
\begin{equation*}
E(t) \ge \fr{1}{2}\int_0^l\left[\left|\fr{\p\u}{\p t}\right|^2 +
\fr{\pi^2 a^2}{l^2}|\u|^2\right]{\rm d}x +
\fr{b}{4}\fr{\pi^4}{l^4}\left[\int_0^l|\u|^2{\rm d}x\right]^2.
\end{equation*}
With the aid of the estimate (6), we obtain
\begin{equation*}
\fr{\pi^2 a^2}{2\,l^2}\int_0^l|\u|^2{\rm d}x +
\fr{b}{4}\fr{\pi^4}{l^4}\left[\int_0^l|\u|^2{\rm d}x\right] \le
ME(0){\rm e}^{-\mu t}
\end{equation*}
whose solution yields the estimate for the total amplitude of the string
\begin{equation}
\int_0^l|\u|^2{\rm d}x \le \fr{l^2 a^2}{\pi^2 b}\left[\sqrt{1 +
\fr{4b}{a^4}ME(0){\rm e}^{-\mu t}} \;- 1\right].
\end{equation}
We say that the estimate (22) is really close for the quasi-steady
stage $t \gg 0$ when the velocities die down. This is because the kinetic
energy part in the derivation of (22) gets ignored. Apparantly when
$t\rightarrow +\,\infty,\;|\u|\rightarrow 0$.\vspace{-.6pc}

\section{Conclusion}

Here we have established uniform stability of a transversely vibrating
string (or a metallic wire) fixed at both ends, governed by vector
Kirchoff type nonlinear integro-differential equation, which takes
into account the string's elasticity and passive viscous damping. It has
been proved directly from the equations of motion that the energy of the
system decays exponentially with time. Study of nonlinear vibrations
assumes significance in analysing slender structural elements capable of
withstanding finite deformations. This paper is motivated by such
considerations.\vspace{-.6pc}

\section*{Appendix}

Here we give a simple proof of Scheefer's inequality (10). Since $\u(0,
t) = {\bf 0} = \u(l, t)$, we can have a Fourier representation
\begin{equation*}
\u(x, t) = \sum_{n=1}^{\infty}{\bf b}_n(t)\sin\fr{n\pi x}{l}.
\end{equation*}
Assuming term-by-term differentiation holds,
\begin{equation*}
\fr{\p\u}{\p x} = \fr{\pi}{l}\sum_{n=1}^{\infty}n{\bf
b}_n(t)\cos\fr{n\pi x}{l}.
\end{equation*}
Applying Parseval's theorem to the two series, we get
\begin{equation*}
\int_0^l\left|\fr{\p\u}{\p x}\right|^2{\rm d}x =
\fr{l}{2}\fr{\pi^2}{l^2}\sum_{n=1}^{\infty}n^2|{\bf b}_n|^2 \ge
\fr{l}{2}\fr{\pi^2}{l^2}\sum_{n=1}^{\infty}|{\bf b}_n|^2 =
\fr{\pi^2}{l^2}\int_0^l|\u|^2{\rm d}x.
\end{equation*}
Hence the inequality.\vspace{-.6pc}

\section*{Acknowledgement}

The authors wish to thank the referee for critical comments which led to
several clarifications in the text.

\end{document}